\documentclass[reqno, 10pt]{amsart}
\usepackage{amssymb,mathrsfs}
\usepackage[usenames, dvipsnames]{color}
\usepackage{esint}
\usepackage[colorlinks=true,linkcolor=blue,citecolor=red,urlcolor=black]{hyperref}
\usepackage{todonotes}
\usepackage{verbatim}
\usepackage{enumerate}
\usepackage{bbm}

\usepackage{cite}
\usepackage[nobysame,non-sorted-cites,initials]{amsrefs}
\def\MR#1{}

\theoremstyle{plain}
\newtheorem{theorem}{Theorem}[section]

\newtheorem{proposition}[theorem]{Proposition}

\theoremstyle{definition}

\theoremstyle{remark}

\newcommand{\dv}{\operatorname{div}}
\newcommand{\osc}{\operatorname{osc}}

\newcommand{\dist}{\operatorname{dist}}

\newcommand{\sgn}{\operatorname{sgn}}

\numberwithin{equation}{section}

\newcommand{\bC}{\mathbb{C}}
\newcommand{\bN}{\mathbb{N}}

\newcommand{\bR}{\mathbb{R}}

\newcommand{\bS}{\mathbb{S}}

\newcommand\cA{\mathcal{A}}
\newcommand\cB{\mathcal{B}}

\newcommand\cD{\mathcal{D}}

\newcommand\cF{\mathcal{F}}

\newcommand\sA{\mathscr{A}}

\makeatletter
\def\dashint{\operatorname%
{\,\,\text{\bf--}\kern-.98em\DOTSI\intop\ilimits@\!\!}}
\makeatother

\newcommand\dashnorm[2]{\nparallel\kern-.2em #1 \Vert_{#2}}

\usetikzlibrary{patterns}

\begin{document}
\title[Recent developments on elliptic equations from composites]{Recent developments on elliptic equations from composites}

\author[H. Dong]{Hongjie Dong}

\author[Z. Yang]{Zhuolun Yang}

\address[H. Dong]{Division of Applied Mathematics, Brown University, 182 George Street, Providence, RI 02912, USA}
\email{Hongjie\_Dong@brown.edu}

\address[Z. Yang]{Department of Mathematics, The Ohio State University, 231 West 18th Avenue, Columbus, OH 43210, USA}
\email{yang.8242@osu.edu}

\thanks{H. Dong was partially supported by the NSF under agreements DMS-2055244 and DMS-2350129.}
\thanks{Z. Yang was partially supported by the NSF under agreements DMS-2550221.}


\date{\today}

\begin{abstract}
When inclusions in a composite are separated by a very small gap, high contrast between the inclusion and matrix properties can induce strong amplification of the underlying field inside the narrow region. Quantifying this field concentration phenomenon is important both for the theory of composite materials and for practical applications. This survey reviews substantial progress over the past three decades. In particular, we survey a set of elliptic equations and systems for which optimal estimates or sharp asymptotic characterizations have been obtained, and we highlight several interesting open questions.
\end{abstract}

\maketitle
\tableofcontents

\section{Introduction}

Composites are prevalent both in nature and in engineered materials. Strengthening fibers are embedded within a matrix that protects them and transfers loads between fibers, resulting in composites with properties superior to those of the individual components. Many models for composites can be reduced to partial differential equations with coefficients that jump across material interfaces. A common and practically important scenario is when two inclusions come very close to one another. Even if the boundary data are smooth, the narrow gap between inclusions can force the solution to change rapidly across a short distance, which may create very large gradient fields. Understanding the field concentration phenomenon quantitatively is important because it can trigger material failure (see, for example, \cites{BASL,Kel}). In some other cases, inclusions are deliberately designed to create the field concentration to achieve desired enhancement of the field (see, for example, \cites{acsphotonics,KANG20191670}). 

Theoretical analysis of field concentration phenomenon traces back to effective medium theory. In the pioneering work \cite{Keller}, the electrical or thermal resistance of the narrow gaps between the inclusions was analyzed in order to estimate the effective conductivity of a medium dense arrays of nearly touching  perfectly conducting or insulating cylinders (2D) or perfectly conducting spheres (3D). See also \cite{Keller2} for related analysis in elasticity. Over the past three decades, substantial progress has been made in quantitatively characterizing field concentration. In this article, we survey some of the important recent progress.

\begin{figure}[h]
\begin{tikzpicture}[scale=1.0]\label{fig:domain}
    \draw[smooth cycle, tension=0.9]
        plot[domain=0:360, samples=200]
        ({20mm*(0.9+0.12*cos(3*\x)+0.07*sin(5*\x))*cos(\x)},
         {20mm*(0.9+0.12*cos(3*\x)+0.07*sin(5*\x))*sin(\x)});
    \draw (0, 0.55) circle [radius=5mm];
    \draw (0, -0.55) circle [radius=5mm];
    \node (core) at (0.5, 0) {$\varepsilon$};
    \node (D_1) at (0, 0.55) {$D_1$};
    \node (D_2) at (0, -0.55) {$D_2$};
    \node (shell) at (30:13mm) {$\Omega$};
\end{tikzpicture}
\caption{Domain $\Omega$}
\end{figure}
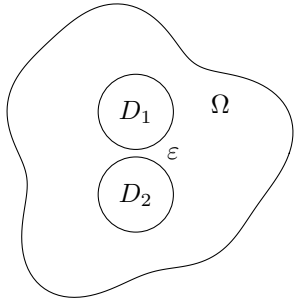

Since the strongest concentrating effect is generated near the local narrow region between two inclusions, it is natural to isolate the pair that forms the narrow neck and study the problem in between. Mathematically, the nature of the domain is as follows: Let $\Omega$ be a bounded open domain in $\bR^n$ with $C^{2}$ boundary, and let $D_{1}$ and $D_{2}$ be two strictly convex open sets whose closure belong to $\Omega$, with $\varepsilon$ distance apart. See figure \ref{fig:domain} for example.

We restrict our attention to elliptic conductivity problem and to the Lam\'e system of elasticity. The following second-order elliptic equations in divergence form with discontinuous coefficients model the conductivity problem when inclusions are \textit{perfectly bonded} into the background matrix:
\begin{equation}\label{equk}
\begin{cases}
\mathrm{div}\Big(a(x)\nabla{u}(x)\Big)=0&\mbox{in}~\Omega,\\
u=\varphi(x)&\mbox{on}~\partial\Omega,
\end{cases}
\end{equation}
where
\begin{equation}\label{conductivity}
a(x) = k_1 \chi_{D_1} + k_2 \chi_{D_2} + \chi_{\Omega\setminus \overline{(D_1 \cup D_2)}}.
\end{equation}
$a(x)$ is viewed as conductivity, the solution $u$ represents the voltage potential, and the associated physical field $E = -\nabla u$ is the electric field. A key feature of the perfectly bonded problem is the continuity of potential and the continuity of flux across the interfaces:
\begin{equation}\label{transmission}
u |_+ = u|_- \quad \mbox{and}\quad \left.\frac{\partial{u}}{\partial\nu} \right|_+ = k_i \left.\frac{\partial{u}}{\partial\nu} \right|_- \quad \mbox{on}~\partial D_i.
\end{equation}
Here and throughout the paper, the subscripts
$\pm$ indicate the limit from outside and inside the inclusion, respectively, and $\nu$ denotes the inner normal vector on $\partial {D}_{1} \cup \partial {D}_{2}$. In section \ref{sec:robin}, we consider a model in which the bonding between the inclusions and the background matrix is imperfect. As a result, one of the transmission conditions \eqref{transmission} does not hold. 

Another important class of models comes from elasticity. In that setting, one studies the Lam\'e system
\begin{equation}\label{lame}
\begin{cases}
\mathrm{div}\Big( (\chi_{\Omega\setminus \overline{(D_1 \cup D_2)}} \bC^0 + \chi_{D_1 \cup D_2} \bC^1) e(u) \Big)=0&\mbox{in}~\Omega,\\
u=\varphi(x)&\mbox{on}~\partial\Omega,
\end{cases}
\end{equation}
where $u = (u^{(1)}, u^{(2)}, \cdots, u^{(n)})^T$ denotes the displacement field, $$e(u) = \frac{1}{2} (\nabla u + (\nabla  u)^T)$$ is the strain tensor. The tensors $\bC^0$ and $\bC^1$ are the elasticity tensors for the background and the inclusion, respectively, and are given by
\begin{equation}\label{C0_tensor}
C^0_{ijkl} = \lambda \delta_{ij} \delta_{kl} + \mu (\delta_{ik} \delta_{jl} + \delta_{il} \delta_{jk}), \quad C^1_{ijkl} = \lambda_1 \delta_{ij} \delta_{kl} + \mu_1 (\delta_{ik} \delta_{jl} + \delta_{il} \delta_{jk}),
\end{equation}
where $(\lambda,\mu)$ and $(\lambda_1, \mu_1)$ are distinct Lam\'e constants, and $\delta_{ij}$ denotes the Kronecker symbol.

In a seminal work, Babu\v{s}ka, Andersson, Smith, and Levin \cite{BASL} studied the Lam\'e system \eqref{lame} and provided numerical evidence that, when the Lam\'e constants $(\lambda,\mu)$ and $(\lambda_1, \mu_1)$ are bounded away from $0$ and infinity, the gradient of solutions remain bounded uniformly in $\varepsilon$. Later, Bonnetier and Vogelius \cite{BV} studied the conductivity problem \eqref{equk}, and proved that for fixed $k_1, k_2$, the electric field $|\nabla u|$ is bounded when $\varepsilon = 0$, $n = 2$ and the inclusions $D_1$ and $D_2$ are disks. This result was extended by Li and Vogelius \cite{LV} to general second order elliptic equations of divergence form with piecewise H\"older coefficients and general shape of inclusions $D_1$ and $D_2$ in any dimension. Furthermore, they established a stronger piecewise $C^{1,\alpha}$ control of $u$, which is independent of $\varepsilon$. Li and Nirenberg \cite{LN} further generalized these global Lipschitz and piecewise $C^{1,\alpha}$ estimates to second-order elliptic systems of divergence form, including the Lam\'e system of elasticity \eqref{lame}. Taken together, these works suggest that strong field concentration arises only in high-contrast regimes, for instance when $k_1$ or $k_2$ tends to $0$ or $\infty$. In such high-contrast cases, numerical studies \cites{Kel, BudCar, Mar} indicate that $|\nabla u|$ typically diverges as $\varepsilon\to 0$.

This survey focuses on several representative high-contrast models for which substantial progress has recently been made. Our goal is not only to collect sharp estimates, but also to highlight several qualitative mechanisms that distinguish the different models. In the linear perfect conductivity problem, the optimal blow-up rate is by now well understood and, remarkably, is largely independent of the local geometry of the inclusions. By contrast, for the insulated problem in dimensions $n\ge 3$, the optimal exponent depends on the local geometry through an eigenvalue problem on the sphere. In nonlinear insulated problems involving the $p$-Laplacian, threshold phenomena appear and the behavior changes depending on the relation between $p$ and the dimension. For imperfect bonding interfaces, the interfacial effect may introduce an additional lower-order term, leading to a new dichotomy between boundedness and blow-up.

The paper is organized as follows. In Section \ref{sec:perfect} we review the perfect conductivity problem, its elasticity analogue with hard inclusions, and the nonlinear perfect conductivity problem for the $p$-Laplacian. Section \ref{sec:insulated} is devoted to the insulated conductivity problem in both linear and nonlinear settings. In Section \ref{sec:mix} we discuss transmission problems with coefficients of different signs. Section \ref{sec:robin} concerns conductivity problems with imperfect bonding interfaces, with particular emphasis on low-conductivity-type interfaces. We conclude in Section \ref{sec:conclude} with some final remarks and a collection of open problems. For a complementary overview of the area, we refer the reader to the survey article \cite{Kang}.

\subsection{Notation}

We use the notation $x = (x', x_n)$, where $x' \in \bR^{n-1}$. After choosing a suitable coordinate system, we may assume that the
shortest segment connecting $\partial D_1$ and $\partial D_2$ is centered at the origin
and is parallel to the \(x_n\)-axis. In these coordinates, the portions of the two
boundaries facing the narrow gap, denoted by $\Gamma_+$ and $\Gamma_-$, can be
represented near the origin as graphs
\begin{align*}
\Gamma_+ = \left\{ x_n = \frac{\varepsilon}{2}+f(x')\right\} \quad \mbox{and} \quad\Gamma_- = \left\{ x_n = -\frac{\varepsilon}{2}+g(x')\right\},
\end{align*}
where $f$ and $g$ are $C^{1,1}$ functions satisfying
\begin{align*}
&f(0')=g(0')=0, \quad \nabla_{x'}f(0')=\nabla_{x'}g(0')=0,\\
&f(x)-g(x) > c_1|x'|^2, \quad \|f\|_{C^{1,1}} + \|g\|_{C^{1,1}} \le c_2,
\end{align*}
with some positive $c_1$ and $c_2$.

\begin{figure}[h]\label{fig:narrow}
\begin{tikzpicture}[x=1cm,y=1cm,>=stealth,scale=0.6]

\def\a{0.55}        
\def\c{0.02}       
\def\xmin{-6.6}
\def\xmax{ 6.6}
\def\xA{ 2.2}       
\def\xB{ 6.1}       

\newcommand{\Top}[1]{\a+\c*(#1)^2}
\newcommand{\Bot}[1]{-\a-\c*(#1)^2}

\fill[
  pattern=north east lines,
  pattern color=cyan!70,
]
  plot[domain=\xA:\xB,samples=120] (\x,{\Top{\x}})
  -- plot[domain=\xB:\xA,samples=120] (\x,{\Bot{\x}})
  -- cycle;

\draw[gray!65,line width=0.8pt] (\xA,{\Top{\xA}}) -- (\xA,{\Bot{\xA}});
\draw[gray!65,line width=1.2pt] (\xB,{\Top{\xB}}) -- (\xB,{\Bot{\xB}});

\draw[line width=1.2pt]
  plot[domain=\xmin:\xmax,samples=220] (\x,{\Top{\x}});
\draw[line width=1.2pt]
  plot[domain=\xmin:\xmax,samples=220] (\x,{\Bot{\x}});

\node[font=\Large] at (-1,1.75)
  {$\Gamma_{+}=\left\{x_n=\frac{\varepsilon}{2}+f(x')\right\}$};
\node[font=\Large] at (-0.7,-1.5)
  {$\Gamma_{-}=\left\{x_n=-\frac{\varepsilon}{2}+g(x')\right\}$};

\coordinate (x0) at (4.15,0.18);
\fill (x0) circle (2.0pt);
\node[font=\Large,below=2pt] at (x0) {$x_0$};

\draw[<->,line width=1.0pt]
  (4.45,0.18) -- node[above,font=\Large] {$r$} (5.85,0.18);

\end{tikzpicture}
\caption{The narrow gap between two inclusions.}
\end{figure}
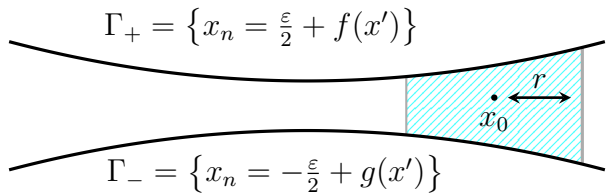
We note that in some models, inclusions with less regular boundary have also been studied in the literature. See, for instance, \cites{KANG20191670, MR4266231}.

We also use the notation
$$
\widetilde{\Omega}:=\Omega \setminus \overline{(D_1 \cup D_2)}
$$
for the background region outside the inclusions. For $0 < r<1$ and $x_0 \in \widetilde\Omega$, we denote
\begin{align*}
\Omega_{x_0,r}:=\left\{(x',x_{n})\in \widetilde{\Omega}~\big|~-\frac{\varepsilon}{2}+g(x')<x_{n}<\frac{\varepsilon}{2}+f(x'),~|x' - x_0'|<r\right\}
\end{align*}
to be the local cylindrical region centered at $x_0$, indicated in Figure \ref{fig:narrow}. We also denote
$$
\Omega_{r}:= \Omega_{0,r}.
$$
This notation will be used repeatedly in the local analysis near the gap.


\section{Perfect conductivity problem, elasticity with hard inclusions}\label{sec:perfect}


\subsection{Linear perfect conductivity problem}

We begin with the classical high-contrast limit in which both inclusions are
perfect conductors, that is,
$$
k_1=k_2=\infty.
$$
In this regime, the potential is constant on each inclusion, and the conductivity
problem reduces to
\begin{equation}\label{perfect}
\begin{cases}
\Delta u =0&\mbox{in}~\widetilde\Omega,\\
u  = U_j &\mbox{on}~\partial D_j, j = 1,2,\\
\int_{\partial D_j} \partial_\nu u \, d\sigma = 0& j=1,2,\\
u=\varphi&\mbox{on}~\partial\Omega,
\end{cases}
\end{equation}
where the constants \(U_1\) and \(U_2\) are determined by the zero-flux conditions. This formulation can be derived from \eqref{equk}--\eqref{transmission}. See, for
instance, the appendix of \cite{BLY1}. There is also a series of work treating the case $\Omega =\bR^n$, in which the boundary condition \eqref{perfect}$_4$ is replaced by the far-field condition
\begin{equation}\label{perfect_newbc}
u - H = O(|x|^{1-n}) \quad \mbox{as}~ |x| \to \infty,
\end{equation}
where $H$ is a prescribed harmonic function.

Among the models surveyed in this paper, the linear perfect conductivity problem
is the most thoroughly understood. In particular, the optimal blow-up rate of
$|\nabla u|$ as the inclusion distance $\varepsilon\to 0$ is now known in all dimensions, and, remarkably, it depends only on the dimension rather than on the detailed local geometry of the inclusions, provided they remain relatively strictly convex.
The first sharp result was obtained in two dimensions by Ammari, Kang, Lee,
Lee, and Lim \cite{AKLLL} and by Ammari, Kang, and Lim \cite{AKL} for two disks of
comparable size in $\Omega=\mathbb{R}^{2}$, where the electric field blows up at the rate $\varepsilon^{-1/2}$. This was later extended by Yun \cites{Y1,Y2} to general bounded strictly convex smooth inclusions in two dimensions. The two-dimensional estimates were later localized and extended to higher dimensions by Bao, Li, and Yin \cites{BLY1,BLY2}. They proved in \cite{BLY1} that
\begin{equation}\label{optimal_perfect}
\begin{cases}
\| \nabla u \|_{L^\infty(\widetilde{\Omega})} \le C\varepsilon^{-1/2} \|\varphi\|_{C^2(\partial \Omega)} &\mbox{when}~n=2,\\
\| \nabla u \|_{L^\infty(\widetilde{\Omega})} \le C|\varepsilon \ln \varepsilon|^{-1} \|\varphi\|_{C^2(\partial \Omega)} &\mbox{when}~n=3,\\
\| \nabla u \|_{L^\infty(\widetilde{\Omega})} \le C\varepsilon^{-1} \|\varphi\|_{C^2(\partial \Omega)} &\mbox{when}~n\ge 4.
\end{cases}
\end{equation}
These estimates were also shown to be optimal.

Since the study of blow-up rates is rather thorough, a detailed characterization of the blow-up profile would be valuable from both mathematical and engineering perspectives. Kang, Lim, and Yun in \cites{KLY1, KLY2} derived asymptotic formulas for two spherical inclusions embedded in $\bR^2$ and $\bR^3$, respectively. Specifically, for two adjacent circular inclusions $D_{1}$ and $D_{2}$ in $\mathbb{R}^{2}$ of radius $r_{1}$ and $r_{2}$ with $\varepsilon$ apart, they proved the following
\begin{theorem}\label{thm:2.1}
Let $R_{j}$, $j=1,2$, be the reflection with respect to $\partial{D}_{j}$, $p_{1}\in{D}_{1}$ and $p_{2}\in{D}_{2}$ be the unique fixed points of $R_{1}R_{2}$ and $R_{2}R_{1}$ respectively, $\vec{n}$ be the unit vector in the direction of $p_{2}-p_{1}$, and $p$ be the middle point of the shortest line segment connecting $\partial{D}_{1}$ and $\partial{D}_{2}$. For a harmonic function $H$ in $\mathbb{R}^{2}$, let $u$ be the solution to \eqref{perfect}$_{1,2,3}$ with $\Omega = \bR^2$ and \eqref{perfect_newbc}. Then
$$\nabla u(x)=\frac{2r_{1}r_{2}}{r_{1}+r_{2}}(\vec{n}\cdot\nabla{H})(\,p)\Big(\frac{x-p_{1}}{|x-p_{1}|^2}-\frac{x-p_{2}}{|x-p_{2}|^2}\Big)+O(1).$$
\end{theorem}
This result was extended by Ammari et al. \cite{ACKLY}, to the case when two inclusions are strictly convex in $\bR^2$. In $\mathbb{R}^{3}$, an analogous characterization is obtained by Kang, Lim, and Yun in \cite{KLY2} in the narrow region between two balls with the same radius $r$ and when $\sqrt{x_{1}^{2}+x_{2}^{2}}\leq r|\log \varepsilon|^{-2}$. The case for different radii was generalized in \cite{LiWangXu}.

While Theorem \ref{thm:2.1} identifies the leading singular term for highly
symmetric geometries, one would like to have an analogous description for more
general inclusions. A finer asymptotic expansion in this direction was obtained
in \cite{LLY}. To state the result, let $\bar{u}\in{C}^{k,1}(\widetilde{\Omega})$ be an auxiliary function satisfying
$$
\bar u=1 \quad \mbox{on}~\partial D_1,
\qquad
\bar u=0 \quad \mbox{on}~\partial D_2\cup\partial\Omega,
$$
and
\begin{align*}
\bar{u}(x)
=\frac{x_{n}-g(x')+\frac{\varepsilon}{2}}{\varepsilon+f(x')-g(x')}\quad\mbox{in}~~\Omega_{R_{0}},
\end{align*}
together with
\begin{equation*}
\|\bar{u}\|_{C^{k,1}(\mathbb{R}^{n}\setminus \Omega_{R_{0}})}\leq\,C.
\end{equation*}
\begin{theorem}
For $\varphi \in C^0(\partial \Omega)$, let $u \in H^1(\widetilde\Omega)$ be the solution to \eqref{perfect}, then
\begin{itemize}
\item[(i)] if $n=2$, $\partial D_{i}$ and $\partial \Omega$ are of $C^{3,1}$, then
\begin{equation*}
\nabla{u}=\frac{Q[\varphi]\sqrt{\varepsilon}}{\Theta}
\nabla\bar{u}+O(1)\|\varphi\|_{C^{0}(\partial\Omega)}\quad\,\mbox{in}~~\widetilde{\Omega};
\end{equation*}
\item[(ii)] if $n=3$, $\partial D_{i}$ and $\partial \Omega$ are of $C^{k,1}$, $k \ge 3$, then
\begin{equation*}
\nabla{u}=\frac{Q[\varphi]}{\Theta}\left(\frac{1}{|\log\varepsilon|-\widetilde{M}}+O(1)\varepsilon^{\frac{1}{2}-\frac{1}{2k}}|\log \varepsilon|^{-1}\right)
\nabla\bar{u}+O(1)\|\varphi\|_{C^{0}(\partial\Omega)}\quad\,\mbox{in}~~\widetilde{\Omega},
\end{equation*}
\end{itemize}
where $Q$ is a specific linear functional on the boundary data $\varphi$, $\Theta$ and $\widetilde{M}$ are some $\varepsilon$-independent constants.
\end{theorem}

One of the main ingredients in \cite{LLY} is
an asymptotic expansion of the Dirichlet energy of the harmonic function $v_{i}$ in $\widetilde{\Omega}$ satisfying
\begin{equation*}
\begin{cases}
\Delta{v}_{i}=0&\mbox{in}~\widetilde{\Omega},\\
v_{i}=\delta_{ij}&\mbox{on}~\partial{D}_{j},~i,j=1,2,\\
v_{i}=0&\mbox{on}~\partial\Omega.
\end{cases}
\end{equation*}
The method in deriving the asymptotics of the gradients are very different from that in \cites{ACKLY,KLY1,KLY2}, which rely more directly on the explicit
structure available for circular inclusions.

Taken together, these results show that the linear perfect conductivity problem is
by now understood not only at the level of optimal blow-up rates, but also, in a
number of important settings, at the level of precise singular asymptotics.


\subsection{Lam\'e system with hard inclusions}

The scalar perfect conductivity problem has a natural elasticity analogue. Instead of letting the scalar conductivity inside the inclusions tend to $+ \infty$, one sends the Lam\'e parameters $(\lambda_1, \mu_1)$ in the inclusions to $+ \infty$. In this limit, the inclusions become rigid: the strain tensor $e(u)$ vanishes in the inclusions, while the displacement remains continuous across the interfaces. In this case, the Lam\'e system \eqref{lame} with \eqref{C0_tensor} becomes
\begin{equation*}\label{lame_hard}
\begin{cases}
\dv \big(\bC^0 e(u)\big) =0&\mbox{in}~\widetilde\Omega,\\
u|_+  = u|_- &\mbox{on}~\partial D_j,~ j = 1,2,\\
e(u) = 0 & \mbox{in}~D_1 \cup D_2,\\
\int_{\partial D_j} \partial_\nu u|_+ \cdot \psi^\alpha \, d\sigma = 0& j=1,2, ~ \alpha = 1,2,\ldots,n(n+1)/2,\\
u=\varphi&\mbox{on}~\partial\Omega.
\end{cases}
\end{equation*}
Here $\{\psi^\alpha\}_{\alpha = 1}^{n(n+1)/2}$ is the basis of the linear space of rigid displacement in $\bR^d$, that is
$$
\{e_i, x_j e_k - x_k e_j\, |\, 1 \le i \le n, 1 \le j < k \le n\},
$$
where $e_1, \ldots, e_n$ denote the standard basis of $\bR^n$.

This problem admits a natural variational characterization. The displacement $u$
is the minimizer of the elastic energy among all admissible deformations that are
rigid inside the inclusions:
$$
I[u] = \min_{v \in \cA} I[v],
$$
where
$$
I[u]:= \frac{1}{2} \int_{\widetilde{\Omega}} \Big( \bC^0 e(v), e(v) \Big) \, dx,
$$
and
$$
\cA = \{ u \in H^1(\Omega; \bR^n): e(u) = 0 ~\mbox{in}~D_1 \cup D_2, u = \varphi~\mbox{on}~\partial \Omega \}.
$$

Compared with the scalar case, the elasticity problem is more delicate because of
its vector-valued structure. Nevertheless, the sharp concentration rates turn out to be the same as in the perfect conductivity problem. Bao, Li, and Li \cites{BLL,BLL2} proved the same upper bounds as in \eqref{optimal_perfect}, and Li \cite{Li21} later established the optimality of these rates. In two dimensions, Kang and Yu \cite{KangYu19} obtained a precise asymptotic
expansion for $\nabla u$ by constructing singular functions from those arising in
the scalar perfect conductivity problem. We refer the reader to Section 3 of \cite{Kang} for a more detailed description of their construction.

Thus, at the level of blow-up rates, the hard-inclusion Lam\'e system closely
parallels the scalar perfect conductivity problem. At the same time, the vectorial
structure introduces additional obstruction, and the available sharp asymptotic results are correspondingly more involved.


\subsection{Perfect conductivity problem with \texorpdfstring{$p$}{}-Laplacian}\label{sec:perfect_p}

We next turn to a nonlinear analogue of the perfect conductivity problem. Here the
inclusions are still assumed to be perfect conductors, but the background matrix
obeys a power-law current-electric field relation
\begin{equation}\label{power_law}
J = \sigma |E|^{p-2} E, \quad p > 1,
\end{equation}
where $J$, $E$, and $\sigma$ represent current, electric field, and conductivity, respectively. This power law is used to model a range of nonlinear media, including certain dielectrics, polymer processing, plasticity, viscous flow models in glaciology, and electro- and thermo-rheological fluids. See \cites{BIK,Idiart,LevKoh,Ruzicka} and the references therein.

The perfect conductivity problem incorporating the power law \eqref{power_law} can be modeled by the following $p$-Laplace equation with $p >1$:
\begin{equation}\label{equinfty}
\left\{
\begin{aligned}
-\dv (|\nabla u|^{p-2} \nabla u) &=0  &&\mbox{in }\widetilde{\Omega},\\
u &= U_j &&\mbox{on}~\overline{D_{j}},~j=1,2,\\
\int_{\partial D_j} |\nabla u|^{p-2} \nabla u \cdot \nu &=0,  && j=1,2,\\
u &= \varphi  &&\mbox{on } \partial \Omega,
\end{aligned}
\right.
\end{equation}
where $U_j$ are constants determined by the third line of \eqref{equinfty}.

Compared to the linear case $p=2$, the nonlinear perfect conductivity problem \eqref{equinfty} was much less understood. For $n\ge 2$, Gorb and Novikov \cite{GorNov} and Ciraolo and Sciammetta \cite{CirSci} proved the following upper bounds:
$$
\|\nabla u\|_{L^\infty(\widetilde \Omega)} \le
\left\{
\begin{aligned}
&C\varepsilon^{-\frac{n-1}{2(p-1)}} && \mbox{when}~p > \frac{n+1}{2},\\
&C\varepsilon^{-1}|\ln \varepsilon|^{\frac{1}{1-p}} && \mbox{when}~p = \frac{n+1}{2},\\
&C\varepsilon^{-1} && \mbox{when}~1<p < \frac{n+1}{2}.
\end{aligned}
\right.
$$
These bounds were shown to be optimal in their respective papers.  Their
arguments rely strongly on the explicit fundamental solution of the $p$-Laplace equation together with the maximum principle.

A more refined question is whether one can go beyond blow-up rates and identify
the leading singular behavior of the solution, in analogy with the linear theory.
This was addressed in a joint work with Zhu \cite{DYZ24}, where the first asymptotic expansion for the
nonlinear perfect conductivity problem was obtained by developing a Schauder
estimate in the narrow region. A key point in that analysis is the structure of the
touching limit $\varepsilon=0$, which reflects the threshold
$p=(n+1)/2$.

Let $D^0_1$ and $D^0_2$ be the translations of $D_1$ and $D_2$ obtained by shifting $D_1$ downward and $D_2$ upward by $\varepsilon/2$, respectively, so that  $D^0_1$ and $D^0_2$ touch at the origin. Set 
$$\widetilde\Omega^0 : = \Omega \setminus \overline{D^0_1 \cup D^0_2}.
$$
Let $u_0$ be the weak $W^{1,p}$-limit of the solutions $u$ to \eqref{equinfty} as $\varepsilon \to 0$. Then the limiting problem changes qualitatively according to whether $p$ is above or below the threshold $(n+1)/2$.

When $p \ge (n+1)/2$, the limiting function $u_0$ takes the same constant
value on the two touching inclusions, but the total flux across $\partial D_1^0$, defined as
\begin{equation}\label{flux}
\cF:= \int_{\partial D_1^0} |\nabla u_0|^{p-2} \nabla u_0 \cdot \nu,
\end{equation}
may be nonzero. More precisely, $u_0$ satisfies
\begin{equation}\label{touching_equation_1}
\left\{
\begin{aligned}
-\dv (|\nabla u_0|^{p-2} \nabla u_0) &=0  &&\mbox{in }\widetilde{\Omega}^0,\\
u_0 &= U_0 &&\mbox{on}~\overline{D_{1}^0}\cup\overline{D_{2}^0},\\
\int_{\partial D_1^0 \cup \partial D_2^0} |\nabla u_0|^{p-2} \nabla u_0 \cdot \nu &=0,\\
 u_0 &= \varphi  &&\mbox{on } \partial \Omega
\end{aligned}
\right.
\end{equation}
for some constant $U_0$.

By contrast, when $p < (n+1)/2$, the limiting function may take different
constant values on the two inclusions, but the flux $\cF$ must vanish. In that case, $u_0$ satisfies
\begin{equation}\label{touching_equation_2}
\left\{
\begin{aligned}
-\dv (|\nabla u_0|^{p-2} \nabla u_0) &=0  &&\mbox{in }\widetilde{\Omega}^0,\\
u_0 &= U_i &&\mbox{on}~\overline{D_{i}^0}\setminus\{0\},~i=1,2,\\
\int_{\partial D_i^0} |\nabla u_0|^{p-2} \nabla u_0 \cdot \nu &=0, &&i=1,2,\\
 u_0 &= \varphi  &&\mbox{on } \partial \Omega
\end{aligned}
\right.
\end{equation}
for constants $U_1$ and $U_2$.

This distinction in the touching limit is one of the key structural features of the
perfect conductivity problem. It underlies the asymptotic expansion in
\cite{DYZ24}, which we now state.

\begin{theorem}\label{thm:expansion}
Let $p>1$, $n \ge 2$, $u \in W^{1,p}(\Omega)$ be the solution of \eqref{equinfty}, $u_0$ be the solution of \eqref{touching_equation_1} or \eqref{touching_equation_2}, $\cF$ be given in \eqref{flux}, $U_1$, $U_2$ be the constants in \eqref{touching_equation_2}$_2$.
We denote
\begin{equation*}\label{Theta}
\Theta(\varepsilon):= \left\{
\begin{aligned}
& \varepsilon^\frac{2p-n-1}{2(p-1)}, && p > \frac{n+1}{2},\\
& |\ln \varepsilon|^{-\frac{1}{p-1}},  && p = \frac{n+1}{2},\\
&1, && 1<p < \frac{n+1}{2}.
\end{aligned}
\right.
\end{equation*} 
Then the following asymptotic expansions hold.
\begin{itemize}
\itemindent=-13pt
    \item[(i)] If $p\ge (n+1)/2$, for $\varepsilon \in(0,1)$ and $x\in \Omega_{1/4}$, we have
\begin{equation*}\label{eq:expansion1}
\nabla u(x)=\Big(0',\frac{\Theta(\varepsilon)}{\varepsilon + f(x') - g(x')}\sgn(\cF)(K|\cF|)^{1/(p-1)}\Big)+ h.o.t,
\end{equation*}
where $K$ is some explicit constant depending only on the principal curvatures of $\partial D_1, \partial D_2$, $n$, and $p$.
\item[(ii)] If $1<p< (n+1)/2$, for $\varepsilon \in(0,1)$ and $x\in \Omega_{1/4}^\varepsilon$, we have
\begin{equation*}
    \nabla u(x)=\Big(0',\frac{U_1-U_2}{\varepsilon + f(x') - g(x')}\Big)+h.o.t.
\end{equation*}
\end{itemize}
\end{theorem}

The significance of Theorem \ref{thm:expansion} is that it identifies the leading
singular term in the nonlinear perfect conductivity problem without relying on the
explicit fundamental solution of the \(p\)-Laplace equation. Instead, the argument
is based on a Schauder theory adapted to the narrow-gap geometry. This provides
a different route from the earlier blow-up estimates in \cites{GorNov,CirSci}, although the proof still uses the maximum principle at one stage. It would be very interesting to develop an alternative approach that avoids the maximum principle and, more importantly, could be extended to nonlinear systems.


\section{Insulated conductivity problem}\label{sec:insulated}


\subsection{Linear insulated conductivity problem}\label{sec:linear_insulated} 

We next consider the opposite extreme regime
\[
k_1=k_2=0,
\]
for which the inclusions are perfect insulators. In this case, the conductivity
problem \eqref{equk} reduces to
\begin{equation}\label{insulated_1}
\left\{
\begin{aligned}
-\Delta u &=0  &&\mbox{in } \widetilde{\Omega},\\
\frac{\partial u}{\partial \nu} &= 0  &&\mbox{on } \partial (D_1 \cup D_2),\\
 u &= \varphi  &&\mbox{on } \partial \Omega,
\end{aligned}
\right.
\end{equation}
which is known as the insulated conductivity problem. This reduction is explained, for instance, in the appendix of \cite{BLY2}.

Compared with the linear perfect conductivity problem, progress for the insulated problem has been slower and less comprehensive. Using a harmonic conjugate argument, Ammari et al. \cites{AKL,AKLLL} showed that in two dimensions the insulated problem has the same blow-up rate, $\varepsilon^{-1/2}$, as the perfect conductivity problem. In all dimensions $n\ge 2$, Bao, Li, and Yin \cite{BLY2} obtained the upper bound
$$
\|\nabla u\|_{L^\infty(\widetilde{\Omega})} \lesssim \varepsilon^{-1/2}.
$$
In three dimensions, Yun \cite{Y3} proved a sharper blow-up rate restricted on the shortest line segment connecting two inclusions: if $D_1, D_2$ are unit balls, then the blow-up rate on this segment is of order 
$$\varepsilon^{-1+\sqrt{2}/2}.$$ 
Progress in dimensions $n\ge 3$ then stalled for several years, until \cite{LY2}, where Li and the second-named author used Harnack’s inequality and the maximum principle to improve the general upper bound to $\varepsilon^{-1/2+\beta}$ for some $\beta>0$. Shortly thereafter, Weinkove \cite{We} used a Bernstein-type argument to produce a more explicit constant $\beta(n)$ for $n \ge 4$ when the inclusions are balls. Meanwhile, identifying the optimal blow-up rate in dimensions $n \ge 3$ was highlighted by Kang in his ICM 2022 lecture \cite{Kang} as one of the central open questions in the area.

The optimal blow-up rate was finally identified in our joint work with Li \cites{DLY,DLY2}. A striking feature of this result is that, unlike the perfect conductivity problem, the optimal exponent for the insulated problem in dimensions $n\ge 3$ depends on the local geometry of the inclusions. More precisely, it is determined by an eigenvalue problem on $\bS^{n-2}$ involving the principal curvatures of the two boundaries at the closest points.

To explain the result, we focus on the local problem near the origin,
\begin{equation}\label{narrow}
\left\{
\begin{aligned}
-\Delta u &=0 \quad \mbox{in }\Omega_{1},\\
\frac{\partial u}{\partial \nu} &= 0 \quad \mbox{on } \Gamma_+ \cup \Gamma_-.
\end{aligned}
\right.
\end{equation}
In \cite{DLY}, we considered the case when $f,g\in C^{2,\gamma}$ and 
$$
f(x')-g(x') =  a |x'|^2 + O(|x'|^{2+\gamma})\quad\mbox{for}~~0<|x'|< 1,
$$
with $a > 0$, $0 < \gamma < 1$. This includes the case when two inclusions are balls. The following optimal estimate was proved in \cite{DLY}.
\begin{theorem}
For $n \ge 3$, $\varepsilon \in (0,1/4)$, let $u \in H^1(\Omega_{1})$ be a solution of \eqref{narrow}. Then 
\begin{equation}\label{optimal_gradient}
|\nabla u(x)|\lesssim_{n,\gamma,a} \|u\|_{L^\infty(\Omega_{1})} (\varepsilon + |x'|^2)^{\frac{\alpha-1}{2}}\quad \forall x\in \Omega_{1/2},
\end{equation}
where
\begin{equation}\label{alpha}
\alpha:= \frac {-(n-1)+\sqrt{(n-1)^2+4(n-2)}}{2}\in (0,1)
\end{equation}
is the optimal exponent.
\end{theorem}
The optimality of \eqref{optimal_gradient} was also established.
\begin{theorem}
\label{optimality}
For $n \ge 3$, $\varepsilon \in (0,1/4)$, let $\Omega = B_5$, and $D_1 , D_2$ be the unit balls centered at $(0', 1 + \varepsilon/2)$ and $(0', -1 - \varepsilon/2)$, respectively. Let $\varphi = x_1$ and $u \in H^1(\widetilde\Omega)$ be the solution of \eqref{insulated_1}. Then 
\begin{equation*}
\label{grad_u_lower_bound}
\|\nabla u\|_{L^\infty(\widetilde\Omega \cap B_{2\sqrt{\varepsilon}})} \gtrsim_n \varepsilon^{\frac{\alpha-1}{2}},
\end{equation*}
where $\alpha$ is given by \eqref{alpha}.
\end{theorem}
The result above already shows that, even for balls, the insulated problem behaves
very differently from the perfect conductivity problem for dimensions $n \ge 3$.
In \cite{DLY2}, this analysis was extended to general strictly convex inclusions. More precisely, assume that
$$
f(x')-g(x') =  \sum_{i=1}^{n-1} a_i x_i^2 + O(|x'|^{2+\gamma})\quad\mbox{for}~~0<|x'|< 1,
$$
with $a_i > 0$, $i = 1, \ldots, n-1$, $0 < \gamma < 1$. Then the optimal
exponent is determined by the eigenvalue problem
\begin{equation}\label{eigenvalue_problem}
-\dv_{\bS^{n-2}}\Big(a(\xi)\nabla_{\bS^{n-2}}u (\xi)\Big) =\lambda a(\xi) u(\xi), \quad \xi \in \bS^{n-2},
\end{equation}
where $a(\xi) = \xi^t D^2(f-g)(0') \xi$. Let $\lambda_1$ be the first nonzero eigenvalue of \eqref{eigenvalue_problem} and
\begin{equation*}\label{alpha2}
\alpha(\lambda_1) := \frac{-(n-1) + \sqrt{(n-1)^2 + 4 \lambda_1}}{2} \in (0,1).    
\end{equation*}
We established gradient estimates similar to \eqref{optimal_gradient} for any $\alpha < \alpha(\lambda_1)$. We also showed that the exponent $\alpha(\lambda_1)$ is optimal, in the sense that estimate \eqref{optimal_gradient} does not hold for $\alpha > \alpha(\lambda_1)$. Finally, the endpoint estimate $\alpha = \alpha(\lambda_1)$ was recently obtained by Li and Zhao \cite{LZ24}.

The main idea in our work can be summarized as follows. Because the Neumann boundary data is zero on the boundaries of both inclusions, the leading order contribution of $\nabla u$ in the narrow gap is expected to be tangential. This motivates us to study the vertical average of the solution $\bar{u}$ in the narrow gap region, i.e.,
$$\bar{u}(x') :=  \fint_{-\varepsilon/2 + g(x') <x_n<\varepsilon/2 + f(x')} u(x',x_n) \, dx_n, \quad |x'| < 1.$$
For the case when inclusions are unit balls, it turns out that $\bar{u}$ satisfies a degenerate equation
\begin{equation*}
\dv((\varepsilon+|x'|^2)\nabla \bar u)=-\sum_{i=1}^{n-1}\partial_i(\overline{a^{in}\partial_n u}) - \sum_{i=1}^{n-1} \partial_i(e^i \partial_i \bar u) \quad\text{in}\,\,B_{1}'\subset \bR^{n-1},
\end{equation*}
where 
$$
|a^{ni}(x)| = |a^{in}(x)| \le C \varepsilon |x'| \quad \mbox{and} \quad |e^i(x')| \le C|x'|^{2+\gamma}.
$$
In order to obtain the optimal estimate of $\nabla \bar u$, we studied this degenerate equation and established the following regularity result. 

\begin{proposition}
For $n \ge 3$, $\sigma > 1$, $\sigma-1 \neq \alpha$,  let $\bar{u} \in H^1(B_{1}')$ be a solution of
\begin{equation*}
\dv\Big[(\varepsilon + r^2)\nabla \bar u\Big]= \dv F + G \quad\text{for}\,\, (r: = |x'|,\xi) \in (0,1)\times \bS^{n-2},
\end{equation*}
where $F$ and $G$ satisfy 
\begin{equation*}
|F(x')| \lesssim r^{\sigma-2}(\varepsilon + r^2) , \quad |G(x')|\lesssim r^{\sigma-3}(\varepsilon + r^2) \quad \mbox{for}~~ x' \in B_1' \subset \bR^{n-1}.
\end{equation*}
Then 
\begin{align*}
\Big(\fint_{\partial B_R'}|\bar u - \bar{u}(0)|^2\,d\sigma\Big)^{1/2} \lesssim R^{\tilde\alpha},
\end{align*}
where $\tilde\alpha := \min \{ \alpha , \sigma - 1\}$, and $\alpha$ is given by \eqref{alpha}.
\end{proposition}
The optimal exponent $\alpha$ comes from the homogeneous part of the degenerate equation. By taking a spherical harmonic decomposition, the first mode $U_1$ of the homogeneous solution satisfies the ODE
$$
U_{1}''(r) +\left(\frac {n-2}r +\frac{2r}{\varepsilon+r^2} \right)U_{1}'(r) - \frac {n-2} {r^2} U_{1}(r)=0, \quad r \in (0,1).  
$$
When $\varepsilon = 0$, the function $r^\alpha$ is an exact solution, and for $\varepsilon>0$, the
maximum principle yields the decay estimate 
$$|U(r)| \le r^\alpha |U(1)|.$$
This is the origin of the optimal exponent in \eqref{optimal_gradient}.

At this point, the linear insulated problem is fairly well understood at the level of sharp blow-up rates. What remains missing is a precise asymptotic expansion for
$\nabla u$, analogous to what is known for the perfect conductivity problem.
Another major challenge is to extend the theory to systems such as the Lam\'e
system: the existing arguments rely strongly on the scalar structure of the equation and on the maximum principle, and these features do not readily persist for the vector-valued system.


\subsection{Insulated conductivity problem with \texorpdfstring{$p$}{}-Laplacian}

We now consider the nonlinear analogue of the insulated conductivity problem. As in Section~\ref{sec:perfect_p}, we assume that the background matrix obeys the power-law current-electric field relation \eqref{power_law}. In that setting, the insulated conductivity problem takes the form
\begin{equation}\label{equzero_p}
\left\{
\begin{aligned}
-\dv (|\nabla u|^{p-2} \nabla u) &=0 \quad \mbox{in }\widetilde{\Omega},\\
\frac{\partial u}{\partial \nu} &= 0 \quad \mbox{on}~\partial {D}_{i},~i=1,2,\\
 u &= \varphi \quad \mbox{on } \partial \Omega.
\end{aligned}
\right.
\end{equation}

Although the nonlinear perfect conductivity problem had already been studied in \cites{GorNov,CirSci}, the nonlinear insulated problem appears to have been considered only recently, in our joint work with Zhu \cite{DYZ23}. Compared with the linear insulated case, the nonlinear theory is much less complete. At present, the most striking feature is a threshold phenomenon: in two dimensions, there is a critical value $p=3$, below which the blow-up rate remains of order $\varepsilon^{-1/2}$, and above which the singularity becomes weaker. The situation in higher dimensions is still not fully understood. 

Our starting point is a universal pointwise estimate, valid for all $p>1$ and all dimensions $n\geq 2$.
\begin{theorem}
\label{thm-1/2}
Let $p>1$, $n \ge 2$, $\varepsilon \in(0,1)$, and $u \in W^{1,p}(\widetilde{\Omega})$ be the solution of \eqref{equzero_p}. Then for any $x\in \Omega_{1/2}$ and $\eta=\frac{1}{4}(\varepsilon+|x'|^2)^{\frac{1}{2}}$,
\begin{equation*}\label{gradient-1/2}
|\nabla u(x)| \lesssim_{n,p} (\varepsilon+|x'|^2)^{-\frac{1}{2}}\underset{\Omega_{x,\eta}}{\osc}~u.
\end{equation*}
\end{theorem}
Thus, exactly as in the linear insulated problem, the natural scale in the narrow
region is $(\varepsilon+|x'|^2)^{\frac{1}{2}}$, rather than the height of the neck $\varepsilon+|x'|^2$. In dimensions $n\geq 3$, we improve the upper bound in Theorem \ref{thm-1/2} to the order of $\varepsilon^{-1/2 + \beta}$ for some $\beta>0$ when $\Omega_{x,\eta}$ is replaced with $\Omega_1$.
\begin{theorem}
\label{Theorem_improved}
Let $p > 1$, $n \ge 3$, $\varepsilon\in(0,1)$, and $u \in W^{1,p}(\widetilde{\Omega})$ be the solution of \eqref{equzero_p}. Then for any $x\in\Omega_{1/2}$,
\begin{equation*}\label{gradient_-1/2+beta}
|\nabla u(x)| \lesssim_{n,p,\beta} (\varepsilon + |x'|^2)^{-1/2 + \beta}\underset{\Omega_{1}}{\osc}~u.
\end{equation*}
\end{theorem}
The estimate above shows that, in dimensions $n\ge 3$, the universal exponent
\(-1/2\) is not optimal. However, the sharp exponent is not yet known in general.
A more explicit result can be obtained when $p > n + 1$, and the inclusion
boundaries satisfy an additional convexity condition. Assume that $f,g\in C^2$ and that
\begin{equation}\label{fg_convex}
\kappa_1I_{n-1} \le \nabla^2 f(x'), -\nabla^2 g(x') \le \kappa_2 I_{n-1} \quad\mbox{for}~~0\leq |x'|<1,
\end{equation}
for some positive constants $\kappa_1$ and $\kappa_2$. Then one has the following
estimate.
\begin{theorem}\label{thm:2d}
Let $n\geq 2$, $p > n+1$, $f$, $g$ be $C^{2}$ functions that further satisfy \eqref{fg_convex}. Let $u \in W^{1,p}(\widetilde{\Omega})$ be the solution of \eqref{equzero_p}. Then for any $\delta>0$ and $\varepsilon\in(0,1)$, we have
\begin{equation*} \label{gradient_2d}
|\nabla u(x)| \lesssim_{n,p,\delta} (\varepsilon + |x'|^2)^{-\frac{n+2\delta}{2(p-1)}}\underset{\Omega_{1}}{\osc}~u \quad \mbox{for}~~x \in \Omega_{1/2}.
\end{equation*}
\end{theorem}

Furthermore,  we show that when $n=2$, the blow-up exponents $-1/2$ for $p\leq 3$ and $-1/(p-1)$ for $p>3$ obtained in Theorems \ref{thm-1/2} and \ref{thm:2d} are critical.
\begin{theorem}\label{thm:2d:2}
For $n =2$, $ p>1$, $\varepsilon\in(0,1)$, let $\Omega = B_5$, and $\cD_1,\cD_2$ be the unit balls centered at $(0,1+\varepsilon/2)$ and $(0,-1-\varepsilon/2)$, respectively. Let $\varphi = x_1$ and $u \in W^{1,p}(\widetilde \Omega)$ be the solution of \eqref{equzero_p}. Then for any $\delta> 0$, when $ p \in(1, 3]$,
$$
\| \nabla u \|_{L^\infty(\widetilde \Omega \cap B_{8\sqrt{\varepsilon/\delta}})} \gtrsim_{p,\delta} \varepsilon^{\frac{-1 + \delta}{2}},
$$
and when $p > 3$,
$$
\| \nabla u \|_{L^\infty(\widetilde \Omega \cap B_{8\sqrt{\varepsilon/\delta}})} \gtrsim_{p,\delta} \varepsilon^{\frac{-1 + \delta}{p-1}}.
$$
\end{theorem}
Finally, we obtain more explicit improvements of the exponent $-1/2$ when the dimension is sufficiently large for any $p>1$.
\begin{theorem}\label{thm:bern}
Let $f$, $g$ be $C^{2,\text{Dini}}$ functions that further satisfy \eqref{fg_convex}. Let $p > 1$, $\beta\in[0,1/2)$, $ \varepsilon \in(0,1)$, and $u \in W^{1,p}(\widetilde \Omega)$ be the solution of \eqref{equzero_p}. If $n$, $p$, and $\beta$ satisfy either
\begin{equation}\label{np_relation_1}
p \ge 2, \quad n \ge \frac{5(p-1)}{2} \left( \frac{p+1 - 2\beta(p-1)}{2} + \frac{ \kappa_2}{(1-2\beta)\kappa_1}  \right) + 1,
\end{equation}
or
\begin{equation}\label{np_relation_2}
1 < p < 2, \quad n \ge \frac{5}{2} \left( \frac{3 - 2\beta}{2} + \frac{ \kappa_2}{(1-2\beta)\kappa_1}  \right) + 3-p,
\end{equation}
then 
\begin{equation}\label{gradient_-1/2_improved}
|\nabla u(x)| \lesssim_{n,\delta,\beta} \|u\|_{L^\infty(\Omega_{1})} (\varepsilon + |x'|^2)^{-\frac{1}{2}+\beta } \quad \mbox{for}~~x \in \Omega_{1/2}.
\end{equation}
\end{theorem}
By \eqref{np_relation_1} and \eqref{np_relation_2}, when $n \to \infty$, $\beta$ can be chosen arbitrarily close to $1/2$. In view of \eqref{gradient_-1/2_improved}, the singularity of $\nabla u$ diminishes as the dimension $n$ increases.

Let us summarize these results in the form of $|\nabla u (0)| \lesssim \varepsilon^{-\alpha}$:
\begin{center}
\begin{tabular}{|l|l|l|}
\hline
Regime & Upper bound exponent &  lower bound exponent \\
\hline
$p\in(1,n+1]$, $n\ge2$ & $\alpha=\frac12$ & $\alpha = \frac12-$ if $n=2$, $p\le3$ \\
$p>n+1$, $n\ge2$& $\alpha = \frac{n}{2(p-1)}+$ & $\alpha = \frac{1}{p-1}-$ if $n=2$, $p>3$ \\
$p > 1$, $n\ge3$& $\alpha=\tfrac12-\beta$ for some $\beta>0$ & \multicolumn{1}{|c|}{---} \\
$p > 1$, $n \to \infty$ & $\alpha=\tfrac12-\beta$ with $\beta\uparrow\tfrac12$ & \multicolumn{1}{|c|}{---} \\
\hline
\end{tabular}
\end{center}

We briefly summarize the main ideas and technical ingredients. The proof of
Theorem \ref{thm-1/2} is based on mean oscillation estimates for $\nabla u$ at scales adapted to the narrow-gap geometry. For Theorem \ref{Theorem_improved}, one uses a special flattening map
in the neck region together with a Krylov-Safonov estimate for the resulting
uniformly elliptic nondivergence-form equation, following the strategy introduced
in \cite{LY2}. For $p>n+1$, we construct barrier functions to prove Theorem~\ref{thm:2d} and the two-dimensional lower bound in Theorem~\ref{thm:2d:2}. The construction uses the $C^2$ regularity and strict convexity of the inclusion boundaries. The high-dimensional estimate in Theorem \ref{thm:bern} is obtained by adapting a Bernstein argument inspired by \cite{We}. A key ingredient there is the classical observation, going back to Uhlenbeck, that $|\nabla u|^q$ is a subsolution of the normalized $p$-Laplace equation for $q \ge p$.

To conclude, the nonlinear insulated problem is currently much less understood
than its linear counterpart. In two dimensions, the threshold at $p=3$ gives a
fairly complete description of the blow-up behavior. However, in dimensions $n\ge 3$, the sharp exponent is unknown. Since the argument for the linear insulated problem relies on a reduction to a degenerate linear equation, it does not appear to carry over to the nonlinear setting. It would be particularly interesting to determine whether $p=n+1$ is the correct threshold in
higher dimensions, in analogy with the two-dimensional threshold $p=3$. 


\section{Conductors with different signs} \label{sec:mix}

Sections \ref{sec:perfect} and \ref{sec:insulated} focus on high-contrast regimes where the gradient field typically concentrates as $\varepsilon=\dist(D_1,D_2) \to 0$. However, there are important mixed configurations in which this
concentration mechanism breaks down. A particularly striking example occurs when
one inclusion approaches the insulating limit while the other approaches the perfect conducting limit. In that case, the gradient may remain bounded even though both phases are individually extreme. This phenomenon was first discovered by Ji and Kang \cite{JiKang}. 

In two dimensions, Ji and Kang studied the conductivity problem \eqref{equk} for
two circular inclusions with
$$
0<k_1<1,
\quad
k_2>1,
$$
using spectral properties of the Neumann-Poincar\'e operator. Among other results,
they proved that
$$
|\nabla^m u| \le C \left( -\frac{(k_1+1)(k_2+1)}{(k_1-1)(k_2-1)}-1 + \sqrt{\varepsilon} \right)^{-m+1}, \quad m\in \bN.
$$
In particular, in the limiting regime $k_1 \to 0$ and $k_2 \to \infty$, this implies
$$
|\nabla^m u| \le C \varepsilon^{-(m-1)/2},\quad m\in \bN.
$$
The estimate shows that the gradient is bounded independent of $\varepsilon$, although higher-order derivatives may still blow up. 

In \cite{DonYan23}, a sharper picture was obtained. We show that, under suitable regularity assumptions on the domain and data, all derivatives $\nabla^m u$ are in fact uniformly bounded for every $m \in \bN$. We briefly describe two cases. 

We first consider the planar case $n = 2$, where $D_1$ and $D_2$ are disks of
radii $r_1$ and $r_2$, respectively. We study the inhomogeneous transmission problem
\begin{equation}\label{eqn:imhomo_trans}
\dv (a(x) \nabla u) = \dv f \quad \mbox{in}~\Omega,
\end{equation}
where $a(x)$ is the piecewise-constant conductivity in \eqref{conductivity}.

\begin{theorem}\label{general_thm}
Let $\varepsilon\in (0,1/2)$, $\mu \in (0,1)$, and $1/2 < r_1, r_2 < 10$. Assume  $0 < k_1 < 1/2$, $k_2 > 2$, and let $u$ be a weak solution of \eqref{eqn:imhomo_trans}. Fix $m \in \bN$ and assume that $f$ is piecewise $C^{2m-1, \mu}$, and for some constant $C_m > 0$,
$$
\|u\|_{L^2(\Omega)} \le C_m, ~~~ \|f\|_{C^{2m-1,\mu}(\widetilde\Omega)} \le C_m, ~~~ \|f\|_{C^{2m-1,\mu}(D_j)} \le C_m \min\{1,k_j\}, j=1,2, 
$$
then
\begin{equation*}
|\nabla^m u(x)| \lesssim_{m,\mu,r_1,r_2} \left\{
\begin{aligned}
&C_m  &&\mbox{in}~~\Omega_{1/2},\\
&\frac{C_m}{k_1+1} &&\mbox{in}~~D_1,\\
&\frac{C_m}{k_2+1} &&\mbox{in}~~D_2.
\end{aligned}
\right.
\end{equation*}
\end{theorem}

To prove Theorem \ref{general_thm}, we use the Green function method developed in \cite{DL}. In contrast to \cite{JiKang}, our argument estimates the $m$-th order derivatives via the $m$-th order finite differences. This allows us to exploit the intrinsic cancellations in the alternating series and, in particular, to avoid an artificial loss in $\varepsilon$. An interesting question is whether one can extend the result in Theorem \ref{general_thm} to higher dimensions $n \ge 3$.

We next turn to general dimensions $n\ge 2$ and to strictly convex inclusions of
general shape. Here we focus on the extreme mixed-contrast configuration
$$
k_1=0,
\qquad
k_2=\infty.
$$
In this limit, the conductivity problem \eqref{equk} reduces to
\begin{equation}\label{eqn:extreme}
\left\{
\begin{aligned}
\Delta{u}=&0 &&\mbox{in}~\widetilde\Omega,\\
\frac{\partial u}{\partial \nu} = &0 && \mbox{on}~\partial{D}_{1},\\
u= &C &&\mbox{on}~\partial{D}_{2},\\
\int_{\partial{\cB}_{2}}\frac{\partial{u}}{\partial\nu}= &0,\\
u = & \varphi&&\mbox{on}~\partial\Omega,
\end{aligned}\right.
\end{equation}
where $C$ is some constant determined by the fourth line of the equation.

In this mixed regime, the behavior is significantly different from both the perfect conductivity and insulated problems discussed earlier. Instead of blow-up in the narrow gap, one obtains uniform control of derivatives of every order.

\begin{theorem}\label{extreme_thm}
Let $u \in H^1(\widetilde\Omega)$ be a weak solution of \eqref{eqn:extreme}. For $m \in \bN$, some constants $\alpha \in (0,1)$, and $C_{m,\alpha} > 0$, if
$$
\|f\|_{C^{m,\alpha}(\{|x'|<1\})} + \|g\|_{C^{m,\alpha}(\{|x'|<1\})} \le C_{m,\alpha},
$$
then there exist constants $\mu \in (0,1)$ and $C$, depending only on $n$, $m$, $\alpha$, and $C_{m,\alpha}$ such that
$$
|\nabla^m u(x)| \le C \mu^{\frac{1}{\sqrt{\varepsilon} + |x'|}} \|u\|_{L^2(\Omega_1)} \quad \mbox{for}~~x \in \Omega_{1/2}.
$$
\end{theorem}

The proof essentially follows an energy argument in \cite{BLLY}, combined with modifications adapted to the mixed boundary conditions. The same conclusion also remains valid when the Laplacian is replaced by the $p$-Laplacian. 

This section illustrates that high contrast alone does not always lead to field concentration. The specific interface configuration matters in an essential way. In the perfect conductivity and insulated problems, the narrow gap amplifies the field, while in the mixed extreme-contrast setting considered here, the interaction between the Dirichlet and Neumann conditions suppresses that mechanism. This provides an important counterpoint to the blow-up phenomena surveyed in the previous sections. 



\section{Conductivity problem with imperfect bonding interfaces}\label{sec:robin}


The preceding sections concern perfectly bonded interfaces, for which the classical transmission conditions \eqref{transmission} hold. However, in many applications, the interfacial effects cannot be neglected. For example, contact resistance created by surface roughness may substantially reduce the effective conductivity of a composite medium. As emphasized in \cite{PRL}, the study of such interfacial effects is challenging both experimentally and theoretically due to the complexity of the underlying microstructure. 

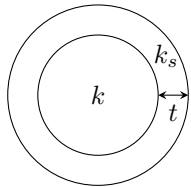
\begin{figure}[ht]
    \centering 
    \begin{tikzpicture}[scale=1.0]
        \draw (0, 0) circle [radius=8mm];
        \draw (0, 0) circle [radius=12mm];
        \node (core) at (0, 0) {$k$};
        \node (shell) at (30:10.5mm) {$k_s$};
        \draw[<->, >=stealth] (8mm, 0)--(12mm, 0);
        \node (width) at (10mm, 0) [below] {$t$};
    \end{tikzpicture}
    \caption{Membrane of core-shell structure}
    \label{core_shell}
\end{figure}

A convenient way to model imperfect bonding is through a thin coating around each
inclusion. consider an inclusion with a core-shell structure, where the shell has a thickness denoted by $t$. Let $k$ and $k_s$ be the conductivities of the core and shell, respectively. See Figure \ref{core_shell}. As $t \to 0$, different scalings of $k_s$ lead to different effective interface laws.

If
$$
\gamma^{-1}: = \lim_{t \to 0} \frac{k_s}{t}
$$
exists and is positive, then the limiting inclusion is said to have an imperfect
bonding interface of low-conductivity type (LC-type), with bonding parameter $\gamma$. In this limit, the transmission conditions on the boundary of the inclusion become
$$
\left.\frac{\partial{u}}{\partial\nu} \right|_+ = k \left.\frac{\partial{u}}{\partial\nu} \right|_- =- \gamma^{-1} (u |_+ - u|_- ) \quad \mbox{on}~\partial \Omega.
$$
On the other hand, if
$$
\alpha = \lim_{t \to 0} k_s t
$$
exists and is positive, then the limiting inclusion is said to have an imperfect
bonding interface of high-conductivity type (HC-type). In that case the transmission conditions take the form
$$
u |_+ = u|_-, \quad  k \left.\frac{\partial{u}}{\partial\nu} \right|_- - \left.\frac{\partial{u}}{\partial\nu} \right|_+ = \alpha \Delta_S u + \nabla_S u \cdot \nabla_S \alpha \quad \mbox{on}~\partial \Omega,
$$
where $\Delta_S$ and $\nabla_S$ are the surface Laplacian and gradient on $\partial \Omega$, respectively. See, for instance, \cites{BenMil}.

In this survey we focus on the case of perfect conductors with LC-type imperfect
bonding interfaces, that is,
$$
k_1=k_2=\infty.
$$
The conductivity problem reduces to the Robin-type boundary value problem
\begin{equation}\label{robin}
\begin{cases}
\Delta u =0&\mbox{in}~\widetilde\Omega,\\
u + \gamma \partial_\nu u = U_j &\mbox{on}~\partial D_j,\, j = 1,2,\\
\int_{\partial D_j} \partial_\nu u \, d\sigma = 0& j=1,2,\\
u=\varphi&\mbox{on}~\partial\Omega,
\end{cases}
\end{equation}
where $\nu$ is the inward normal vector on $\partial D_j$, $U_j$ is some constant determined by the third line of \eqref{robin}. The solution $u\in H^{1}
(\widetilde\Omega)$ to equation \eqref{robin} is indeed the minimizer of a functional in an appropriate function space: $I[u] = \min_{v \in \sA}I[v]$, where
\begin{equation*}
\begin{aligned}
I[v] &:= \int_{\widetilde\Omega} |\nabla v|^2+\gamma^{-1}\int_{\partial D_1} |v-(v)_{\partial D_1}|^2+\gamma^{-1}\int_{\partial D_2} |v-(v)_{\partial D_2}|^2,\\
\sA &:= \{ v \in H^{1}(\widetilde\Omega):\, v = \varphi~~\mbox{on}~~\partial\Omega \},\quad (v)_{\partial D_j}:= \fint_{\partial D_j} v\, d\sigma,\ j=1,2.
\end{aligned}
\end{equation*}

This model was studied in two dimensions by Fukushima et al. \cite{fukushima2024finiteness}, who considered equal disks in $\bR^2$ with far field asymptotic $\varphi=x_2$. They proved that the gradient remains bounded independently of $\varepsilon$, in sharp contrast with the perfect conductivity problem without the membrane, where blow-up occurs as $\varepsilon\to 0$. Motivated by biological membrane models, they conjectured that such boundedness should persist for all $\gamma>0$ and all boundary data.

In \cite{DongYangZhu26}, a recent joint work with Zhu, we showed that the conjecture holds only when the bonding parameter $\gamma$ is sufficiently small, and there is an unexpected dichotomy of field concentration phenomenon depending on the size of $\gamma$. In summary, \cite{DongYangZhu26} establishes three main results:

\begin{enumerate}[1.]
\item \underline{Upper bound for the gradient}: We establish a key upper bound of order $\varepsilon^{-1/2}$ for the gradient of the solution to \eqref{robin}. A direct consequence of this result is that if $\widetilde\Omega$ is symmetric with respect to $x_n$, and the boundary value $\varphi$ is odd in $x_n$, then the gradient is bounded independent of $\varepsilon$. This extends the boundedness result in \cite{fukushima2024finiteness} to all dimensions $n \ge 2$, to more general inclusions $D_1, D_2$, and to more general boundary data $\varphi$.
\item \underline{Absence of field concentration for small $\gamma$}: We show that when the bonding parameter $\gamma$ is sufficiently small, the gradient of the solution is bounded independent of $\varepsilon$. This confirms the conjecture raised in \cite{fukushima2024finiteness} under the additional assumption that $\gamma$ is sufficiently small.
\item \underline{Dichotomy for field concentration}: We show that when $\gamma$ is large, the conjecture fails, and the first boundedness result is highly unstable. More precisely, if the boundary data is perturbed slightly away from being odd in  $x_n$ (for example, $\varphi = \delta x_1 + x_n$ for some small $\delta$), then the gradient of the solution may blow up. In particular, when $D_1, D_2$ are unit balls, $\Omega=B_5$, and $\varphi = x_1$, we prove a dichotomy for the field concentration phenomenon: the gradient of the solution to \eqref{robin} is bounded when $0 < \gamma \le 1$, but blows up when $\gamma > 1$. Moreover, we establish the optimal gradient estimate when $\gamma > 1$. Our analysis shows that in this regime the solution behaves much like the insulated conductivity problem. In fact, as $\gamma \to \infty$, our estimates recover the corresponding results for the perfectly bonded insulated conductivity problem.
\end{enumerate}

The first main result is an anisotropic pointwise gradient estimate in the narrow neck region. As in Theorem \ref{thm-1/2}, we show that the relevant scale is $(\varepsilon+|x'|^2)^{\frac{1}{2}}$, rather than the neck height $\varepsilon+|x'|^2$. This local estimate is fundamental to the analysis throughout the paper \cite{DongYangZhu26}.

\begin{theorem}
\label{thm2-1/2}
Let $n \ge 2$, $\gamma > 0$, $\varepsilon < 1/4$, and $u \in H^1(\widetilde\Omega)$ be the solution of \eqref{robin}. Then, for any $x\in \Omega_{1/2}$ and $r=\frac{1}{4}(\varepsilon+|x'|^2)^{\frac{1}{2}}$, we have
\begin{equation}\label{gradient2-1/2}
|\nabla u(x)| \lesssim_{n,\gamma}  r^{-1} \|u-(U_1 + U_2)/2\|_{L^2_{avg}(\Omega_{x,r})}  + |U_1 - U_2| .
\end{equation}
\end{theorem}
Here
$$
\|f\|_{L^2_{avg} (\Omega)} := \Big( \fint_\Omega |f|^2 \Big)^{1/2}.
$$
A quick consequence is that the gradient remains bounded under the odd symmetry assumption mentioned above. In this case, the symmetry ensures that the constants on the two inclusions cancel, so the singular term in the estimate above vanishes.

The next theorem shows that such boundedness does not require symmetry when the bonding parameter is sufficiently small.

\begin{theorem}\label{thm_polynomial_upperbound}
For any $\beta  \ge 0$,  if $0 < \gamma \ll 1$, then
\begin{equation}\label{grad_u_polyupperbound}
|\nabla u(x)| \lesssim_{n,\gamma,\beta} (\sqrt\varepsilon + |x'|)^\beta \|u - (U_1 + U_2)/2\|_{L^\infty(\Omega_{1})} + |U_1 - U_2| \quad \mbox{in}~ \Omega_{1/2}.
\end{equation}
\end{theorem}

The most interesting phenomenon occurs in the large-$\gamma$ regime. Assume now that $D_1, D_2$ are $C^{2,\sigma}$, and that the gap has the expansion:
$$
f(x') - g(x') = \mu \, |x'|^2 + O(|x'|^{2+\sigma})\quad\mbox{for}~~0<|x'|<1,~ \mu > 0.
$$
Suppose also that, for some horizontal direction $x_j$, $1 \le j \le n-1$, the domain $\widetilde\Omega$ is symmetric with respect to $x_j$ and boundary data is odd in $x_j$. Then the solution is odd in $x_j$, which implies that $U_1 = U_2 = 0$.
Define
\begin{equation}
\label{alpha_n}
\alpha=\alpha(n,\gamma, \mu):= \frac {-(n-1)+\sqrt{(n-1)^2+4(n-2 +  2/(\mu\gamma))}}{2}.
\end{equation}
Then $\alpha <1$ if and only if $\gamma >  1/\mu$. This is exactly the threshold for gradient blow-up.
\begin{theorem}\label{thm_upperbound}
Under the assumption above, we have
\begin{itemize}
\item when $0 < \gamma \le 1/\mu$,
\begin{equation}\label{grad_u_upperbound1}
|\nabla u(x)| \lesssim_{n,\gamma,\mu} \|u\|_{L^\infty(\Omega_{1})} \quad \mbox{in}~ \Omega_{1/2},
\end{equation}
\item when $\gamma >1/\mu$,
\begin{equation}\label{grad_u_upperbound2}
|\nabla u(x)| \lesssim_{n,\mu} \|u\|_{L^\infty(\Omega_{1})} (\varepsilon + |x'|^2)^{\frac{\alpha-1}{2}} \quad \mbox{in}~  \Omega_{1/2},
\end{equation}
\end{itemize}
where $\alpha$ is given in \eqref{alpha_n}.
\end{theorem}

It is further shown that the upper bound in the blow-up regime is sharp.

\begin{theorem}
\label{main_thm}
Let $\Omega = B_{5/\mu}$, $D_1 , D_2$ be balls of radius $1/\mu$ centered at $(0', 1/\mu + \varepsilon/2)$ and $(0', -1/\mu - \varepsilon/2)$, respectively, $\varphi(x) = x_1$, and let $u \in H^1(\widetilde\Omega)$ be the solution of \eqref{robin}. Then there exist positive constants $c$ and $C$, depending only on $n$, $\gamma$, and $\mu$, such that
\begin{equation*}
\|\nabla u\|_{L^\infty(\Omega_{c \sqrt\varepsilon})} \ge \frac{1}{C}\varepsilon^{\frac{\alpha-1}{2}}.
\end{equation*}
\end{theorem}

The threshold also admits a natural interpretation. When $D_1$ and $D_2$ are balls of radius $R$, we have $\mu = 1/R$, so the critical value becomes $\gamma = R$. This is precisely the value for which the inclusions are neutral to uniform fields: when $\gamma = R$, the linear functions $x_j$ automatically satisfy both the Robin interface condition and the zero-flux constraint.

The work in \cite{DongYangZhu26} establishes the first sharp blow-up theory for conductivity problems with imperfect bonding interfaces and reveals a new phenomenon beyond the perfectly bonded models discussed in Sections \ref{sec:perfect}-\ref{sec:mix}: the occurrence of a threshold in the bonding parameter. For small $\gamma$, the membrane suppresses field concentration; for large $\gamma$, the solution behaves more like that of the insulated problem and blow-up may occur.

The proof relies on three main ingredients. The first is the anisotropic gradient estimate \eqref{gradient2-1/2}. To establish it, we construct a flattening map that preserves the Robin boundary conditions on both boundaries, together with a carefully chosen auxiliary function to handle the inhomogeneous boundary condition. This approach is quite robust and applies to a broad class of equations and boundary conditions. Estimate \eqref{grad_u_polyupperbound} is then obtained by an iteration argument, in which the smallness of $\gamma$ yields a suitable decay of $u$ near the origin. Finally, the mechanism behind Theorems \ref{thm_upperbound} and \ref{main_thm} is closely related to that of the insulated problem. After a dimension reduction argument, one reduces the problem to a degenerate elliptic equation in $\mathbb{R}^{n-1}$. The crucial new feature is that the Robin condition contributes an additional zeroth-order term in the reduced equation. It is precisely this term that produces the dichotomy at $\gamma=1/\mu$. In particular, as $\gamma \to \infty$, estimate \eqref{grad_u_upperbound2} recovers the corresponding estimate for the perfectly bonded insulated problem in Section \ref{sec:linear_insulated}, as the constant does not depend on $\gamma$.

On the other hand, as $\gamma \to 0$, the problem \eqref{robin} converges to the perfect conductivity problem \eqref{perfect}. However, estimate \eqref{grad_u_upperbound1} does not carry over as it depends on $\gamma$. The first-named author, together with Li and Zhao \cite{DongLiZhao25}, studied the regime $\gamma \ll 1$ and derived a new estimate that bridges problems \eqref{robin} and \eqref{perfect}.

\begin{theorem}
Let $n \ge 2$, $0 < \gamma \ll 1$, $\varepsilon < 1/4$, for $x \in \Omega_{1/2}$, we have
\begin{equation*}
|\nabla u(x)|\lesssim_n
\left\{ \begin{aligned}
&\frac{1}{\sqrt{\gamma+\varepsilon+|x'|^2}}\|\varphi\|_{C^{2}(\partial\Omega)}& \mbox{for}~n=2,\\
&\frac{1}{ {(\gamma+\varepsilon+|x'|^2)}|\ln(\varepsilon+\gamma)|}\|\varphi\|_{C^{2}(\partial\Omega)}& \mbox{for}~n=3,\\
&\frac{1}{{\gamma+\varepsilon+|x'|^2}}\|\varphi\|_{C^{2}(\partial\Omega)}& \mbox{for}~n\geq 4.
\end{aligned} \right.
\end{equation*}
\end{theorem}
The proof is based on new gradient estimates for elliptic equations, in both $\bR^n$ and the reduced dimension space $\bR^{n-1}$, as $\gamma\to 0$, the maximum principle, and an energy iteration argument in the spirit of \cite{BLLY}. These upper bounds were shown to be optimal in \cite{DongLiZhao25}.

A natural next step is to study inclusions with finite conductivities $k_1$ and $k_2$ with an LC-type imperfect bonding interface. It is well known that, without an LC-type interface, the gradient remains uniformly bounded in $\varepsilon$. See \cites{LV,LN}. It would be interesting to determine whether the addition of an imperfect bonding interface can cause blow-up and, if so, what the optimal blow-up rate is.

By contrast, little is known about HC-type interfaces. In two dimensions, \cite{fukushima2024finiteness} shows that insulators with an HC-type interface and perfect conductors with an LC-type interface are dual to each other. It would be interesting to study HC-type interfaces in dimensions $n \ge 3$, although this appears to be challenging due to the complexity of the transmission conditions.


\section{Concluding remarks and open problems}\label{sec:conclude}

The results surveyed in this article illustrate that field concentration in composite media is governed by a subtle interaction among contrast, geometry, nonlinearity, and interface effects. In the linear perfect conductivity problem, the sharp blow-up rate is by now well understood and is largely independent of the detailed local geometry of the inclusions. By contrast, for the insulated problem in dimensions $n\ge 3$, the optimal exponent depends on the local quadratic geometry of the boundaries through an eigenvalue problem on $\bS^{n-2}$. In the nonlinear insulated problem, threshold phenomena appear, and in the imperfect bonding problem, the Robin interface produces a new dichotomy in the bonding parameter. These features show that, even within the narrow-gap setting, different high-contrast regimes can exhibit significantly different concentration mechanisms.

Although substantial progress has been made, several important questions remain
open.

First, for the linear insulated conductivity problem, the optimal blow-up rate is now known, but a precise asymptotic expansion for the gradient is still missing.

Second, for the insulated problem with $p$-Laplacian, the picture is relatively clear only in two dimensions. In particular,  the optimal blow-up exponent is still unknown for $n\ge 3$. It would be especially interesting to determine whether $p=n+1$ is the correct threshold in higher dimensions.

Third, most of the sharp results discussed here concern scalar equations. Extending the insulated theory to systems, especially to the Lam\'e system of elasticity, remains a major challenge. The existing arguments for the scalar insulated problem rely heavily on the maximum principle and on a reduction to a scalar degenerate elliptic equation, and these tools do not readily extend to vector-valued systems.

Fourth, Section \ref{sec:mix} shows that in certain mixed extreme-contrast configurations the gradient may remain uniformly bounded, in sharp contrast with the blow-up phenomena in the perfect conductivity and insulated problems. While Theorem \ref{general_thm} establishes this behavior in two dimensions, an interesting open problem is whether the same result can be extended to higher dimensions $n\ge 3$. 

Fifth, for imperfect bonding interfaces, the LC-type problem for perfect
conductors is now fairly well understood, but very little is known when the inclusions have finite conductivities $k_1$ and $k_2$. Since the corresponding
perfectly bonded bounded-contrast problem does not exhibit blow-up, it would be
very interesting to understand whether the addition of an imperfect interface can
create field concentration and, if so, to determine the sharp rate.

Finally, HC-type imperfect interfaces remain largely unexplored, especially in
dimensions $n\ge 3$. Even the correct qualitative picture is not yet clear in this setting, partly because the effective transmission condition is substantially more complicated than that in the LC-type case.



\end{document}